\documentclass[12pt]{article}
\usepackage{setspace}
\usepackage{fullpage}
\usepackage[left=1.25in,right=1.25in,top=1in,bottom=1.2in]{geometry}
\usepackage{amsmath, amsfonts, pstricks, multido, pst-node}

\newcommand{\inv}{\mathop{inv}}
\newcommand{\rank}{\mathop{rank}}

 \newtheorem{conj}{Conjecture}[section]

\begin{document}

\title{Formulas for the dimensions of some affine Deligne-Lusztig
Varieties}

\author{Daniel C. Reuman}

\maketitle

\begin{abstract}
Rapoport and Kottwitz defined the affine Deligne-Lusztig varieties
$X_{\tilde{w}}^P(b\sigma)$ of a quasisplit connected reductive
group $G$ over $F = \mathbb{F}_q((t))$ for a parahoric subgroup
$P$. They asked which pairs $(b, \tilde{w})$ give non-empty
varieties, and in these cases what dimensions do these varieties
have. This paper answers these questions for $P=I$ an Iwahori
subgroup, in the cases $b=1$, $G=SL_2$, $SL_3$, $Sp_4$. This
information is used to get a formula for the dimensions of the
$X_{\tilde{w}}^K(\sigma)$  (all shown to be non-empty by Rapoport
and Kottwitz) for the above $G$ that supports a general conjecture
of Rapoport. Here $K$ is a special maximal compact subgroup.
\end{abstract}

\section{Introduction}\label{Intro}

Let $F$ be $\mathbb{F}_q((t))$ with ring of integers
$\mathcal{O}_F$, and let $G$ be a split connected reductive group
over $F$. Let $L$ be the completion of the maximal unramified
extension of $F$, $\bar{\mathbb{F}}_q((t))$.  Let $\sigma$ be the
Frobenius automorphism of $L$ over $F$. Let $\mathcal{B}_n$ be the
affine building for $G(E)$ where $E/F$ is the unramified extension
of degree $n$ in $L$ (so $E=\mathbb{F}_{q^n}((t))$), and let
$\mathcal{B}_{\infty}$ be the affine building for $G(L)$.  Let $T$
be a split torus in $G$, let $B=UT$ be a Borel subgroup, and let
$I$ be an Iwahori in $G(L)$ containing $T(\mathcal{O}_L)$, where
$\mathcal{O}_L$ is the ring of integers of $L$.  Let $A_M$ and
$C_M$ be the correspondingly specified apartment and alcove, which
we assume are in $\mathcal{B}_1$. We will call these the {\em main
apartment} and the {\em main alcove,} respectively.  We assume
that $C_M$ is in the positive Weyl chamber in $A_M$ specified by
$B$.  Let $P \supseteq I$ be a parahoric subgroup of $G(L)$.  If
$b \in G(L)$ then the $\sigma$-conjugacy class of $b$ is
$\{x^{-1}b\sigma (x) : x \in G(L) \}$.  Let $\tilde{W} =
N(L)/T(\mathcal{O}_L)$ be the extended affine Weyl group, and let
$\tilde{W}_P = N(L) \cap P / T(\mathcal{O}_L)$.  Here $N$ is the
normalizer of $T$.

If $\tilde{w} \in \tilde{W}$, then we define, after Rapoport
\cite{R1} and Kottwitz, the {\em generalized affine
Deligne-Lusztig variety $X_{\tilde{w}}^P(b\sigma) = \{x \in G(L)/P
: \inv_P(x, b\sigma(x)) = \tilde{w} \}$}. Here $\inv_P : G(L)/P
\times G(L)/P \longrightarrow P \backslash G(L) / P = \tilde{W}_P
\backslash \tilde{W} / \tilde{W}_P$ is the relative position map
associated to $P$.  Rapoport asked the question of which pairs
$(b, \tilde{w})$ give rise to non-empty sets, and for these pairs,
what is $\dim (X_{\tilde{w}}^P(b\sigma))$ \cite{R1}. Kottwitz and
Rapoport answered the emptiness/non-emptiness part of this
question for $P=K$, the maximal bounded subgroup of $G(L)$
associated to some special vertex $v_M$ of $C_M$ \cite{KR1}.

In Section~\ref{SL3} of this paper we consider the case $G=SL_3$,
$b=1$, $P=I$.  Complete results on emptiness/non-emptiness and
dimension are shown for this case in Figure~\ref{SL3result}. In
Section~\ref{Sp4} we consider $G=Sp_4$, $b=1$, $P=I$.
Emptiness/non-emptiness results and dimension results are in
Figure~\ref{Sp4result}.  The case $G = SL_2$, $b=1$, $P=I$ can be
done using an even simpler version of the same methods.

Kottwitz and Rapoport showed in \cite{KR1} that for general $G$,
$X_{\tilde{w}}^K(\sigma)$ is non-empty for any $\tilde{w}$
corresponding to a dominant cocharacter in the coroot lattice.
Rapoport conjectured a specific formula for the dimension of the
$X_{\tilde{w}}^K(\sigma)$ in \cite{R2}. The knowledge of the
$X_{\tilde{w}}^I(\sigma)$ mentioned in the previous paragraph
gives knowledge of the $X_{\tilde{w}}^P(\sigma)$, so the
dimensions of the $X_{\tilde{w}}^K(\sigma)$ are computed in
Section~\ref{dimK} for $SL_2$, $SL_3$, and $Sp_4$.  The result is
that $\dim(X_{\tilde{w}}^K(\sigma))=\langle \mu , \rho \rangle$,
where $\mu \in X_{*}(T)$ dominant corresponds to $\tilde{w} \in
\tilde{W}_K \backslash \tilde{W} / \tilde{W}_K$, and $\rho$ is
half the sum of the positive roots for $G$.  This supports the
conjecture of Rapoport in \cite{R2}. Preliminary work toward a
proof of this conjecture has been done with Kottwitz.

In Section~\ref{formulasection} we present a formula that
encapsulates part of the results pictured in
Figures~\ref{SL3result} and~\ref{Sp4result}.  The formula also
holds for $SL_2$.  It is too soon to conjecture that this formula
holds for general $G$.  Some results on emptiness/non-emptiness
for $b \neq 1$, $G = SL_2$, $SL_3$, $Sp_4$ can be found in
\cite{Re1}. Section~\ref{RelatedResults} contains a summary of
these results.

This work has significance to the study of the reduction modulo
$p$ of Shimura varieties.  Interested readers should see the
survey article \cite{R2} by Rapoport.

\section{General Methodology}\label{GM}

For this and the next two sections we let $P=I$, so $\tilde{w} \in
\tilde{W}_P \backslash \tilde{W} / \tilde{W}_P = \tilde{W}$, and
we let $X_{\tilde{w}}(\sigma) = X_{\tilde{w}}^I(\sigma)$. For this
section, the group $G$ is simply-connected, so that $I$ is the
stabilizer of $C_M$. First note that if $\tilde{w}C_M \cap C_M$ is
non-empty, then $X_{\tilde{w}}(\sigma)$ can be identified with a
disjoint union of (non-affine) Deligne-Lusztig varieties, whose
structure and dimension are already known \cite{DL1}. Let $v_1$ be
a vertex in $A_M$ and let $v_2$ be a vertex in $\mathcal{B}_1$ in
the same $G(F)$ orbit as $v_1$. We require that $v_1 \notin C_M$.
Let $Q_1$ be the last alcove in a minimal gallery from $C_M$ to
$v_1$, and let $\mathcal{Q}_2$ be the set of all alcoves $Q_2$
containing $v_2$ such that $Q_2$ and $\sigma (Q_2)$ have some
fixed relative position, $p_r$. Note that $Q_1$ does not depend on
the choice of minimal gallery from $C_M$ to $v_1$. We require that
$p_r$ be such that $Q_2 \cap \mathcal{B}_1 = \{ v_2 \}$. We define
{\em the $(v_1 , v_2 , p_r)$-piece of $X_{\tilde{w}}(\sigma)$}
(which may be empty) to be all alcoves $D \subset
\mathcal{B}_{\infty}$ such that there exists $y \in G(L)$ with
$yC_M = D$, $yQ_1 = Q_2$ for some $Q_2 \in \mathcal{Q}_2$ (so
$yv_1 = v_2$), and $\inv(D,\sigma(D))=\tilde{w}$.  The dimension
of $X_{\tilde{w}}(\sigma)$ is the supremum of the dimensions of
its pieces (this supremum could be infinite, {\em a priori},
although we will show that it is finite in the cases we will
consider). We define {\em the $(v_1,v_2,p_r)$-superpiece } to be
the collection of all alcoves $D \subset \mathcal{B}_{\infty}$
such that there exists $y \in G(L)$ with $yC_M = D$ and $yQ_1 =
Q_2$ for some $Q_2 \in \mathcal{Q}_2$.  So the
$(v_1,v_2,p_r)$-superpiece is the disjoint union, over all
$\tilde{w} \in \tilde{W}$, of the $(v_1,v_2,p_r)$-pieces of the
$X_{\tilde{w}}(\sigma)$ (many of which will be empty). The
approach outlined in this paragraph was suggested by Kottwitz, and
is similar to that used in \cite{K3}.

Note that the structure of the $(v_1,v_2,p_r)$-superpiece does not
depend on $v_2$, as long as $v_2$ is some vertex in
$\mathcal{B}_1$ in the same $G(F)$ orbit as $v_1$.  So for each
$(v_1 , p_r)$-pair, we fix an arbitrary vertex $v_2 \in A_M$ in
the same $G(F)$-orbit as $v_1$, and we compute the possible values
of $\inv(D,\sigma(D))$ for $D$ an alcove in $\mathcal{B}_{\infty}$
for which some $y \in G(L)$ gives $yC_M = D$ and $yQ_1= Q_2$ for
some $Q_2 \in \mathcal{Q}_2$. We will discuss how this computation
is carried out for $SL_3$ and $Sp_4$. The results tell us for
which $\tilde{w}$ the $(v_1,v_2,p_r)$-piece of
$X_{\tilde{w}}(\sigma)$ is non-empty. We will also demonstrate a
way of calculating the dimension of each non-empty piece in the
$(v_1,v_2,p_r)$-superpiece, again only for $SL_3$ and $Sp_4$.
Everything we will do also applies to $SL_2$. Aggregating all this
information over all $(v_1,p_r)$-pairs will tell us, for each
piece of each $X_{\tilde{w}}(\sigma)$, whether it is empty or
non-empty, and what its dimension is.  This gives the
emptiness/non-emptiness and dimension of the
$X_{\tilde{w}}(\sigma)$ themselves.

\section{$SL_3$}\label{SL3}

\begin{figure}
\centerline{\psset{xunit=0.233333333333333 in}
\psset{yunit=0.404145188432738 in}

\begin{pspicture*}(0,0.25)(7.2,4.25)


  \multido{\rxb=3+1,\rxt=0.5+0.5,\ryt=2.5+0.5}{5}{ \psline[linewidth=0.05pt](\rxb,0)(\rxt,\ryt) }
  \multido{\rxb=8+1,\rxt=3.5+1}{12}{ \psline[linewidth=0.05pt](\rxb,0)(\rxt,4.5) }
  \multido{\rxb=1+0.5,\rxt=3.5+1,\ryb=2+-0.5}{5}{ \psline[linewidth=0.05pt](\rxb,\ryb)(\rxt,4.5) }
  \multido{\rxb=4+1,\rxt=8.5+1}{11}{ \psline[linewidth=0.05pt](\rxb,0)(\rxt,4.5) }
  \multido{\ny=0+0.5,\ngpa=3+-0.5}{4}{ \psline[linewidth=0.05pt](\ngpa,\ny)(15,\ny) }
  \multido{\ny=2+0.5,\ngpa=0+0.50}{6}{ \psline[linewidth=0.05pt](\ngpa,\ny)(15,\ny) }

\psset{linestyle=dashed} \psline[linewidth=0.05pt](0,2)(2.5,4.5)
\rput[bl](0.25,2){\rule[-2pt]{0pt}{5pt}$C_M$}


\end{pspicture*}}
\caption{The region containing all
vertices $v_1$ that must be considered for $SL_3$}
\label{v1regionSL3}
\end{figure}
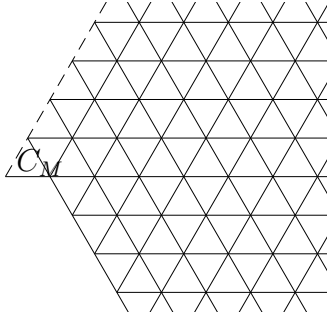

In order to carry out the process outlined in Section~\ref{GM} for
$SL_3$, it suffices to consider $v_1$ in the region pictured in
Figure~\ref{v1regionSL3}.  All other $v_1$ can be obtained from
these by rotating by $120^{\circ}$ or $240^{\circ}$ about the
center of $C_M$.  Further, if $v_1'$ is the rotation of $v_1$ by
$\alpha$ ($=120^{\circ}$ or $240^{\circ}$) about the center of
$C_M$ then it is easy to see that the set $\{ \tilde{w} \in
\tilde{W} : \textrm{the } (v_1',v_2',p_r) \textrm{-piece of }
X_{\tilde{w}}(\sigma) \textrm{ is non-empty} \}$ is the rotation
of the set $\{ \tilde{w} \in \tilde{W} : \textrm{the }
(v_1,v_2,p_r) \textrm{-piece of } X_{\tilde{w}}(\sigma) \textrm{
is non-empty} \}$ by $\alpha$ about the center of $C_M$. Further,
the correspondence between these two sets given by rotation by
$\alpha$ preserves the dimension of the corresponding pieces.

Given the restriction (mentioned in Section~\ref{GM}) that $p_r$
be such that $Q_2 \cap \mathcal{B}_1 = \{ v_2 \}$, we know that
$Q_2$ and $\sigma (Q_2)$ must share exactly one vertex.  So for
$SL_3$, $p_r$ corresponds to some element of $W$, the finite Weyl
group, of length $2$ or $3$.

\begin{figure}
\centerline{\input{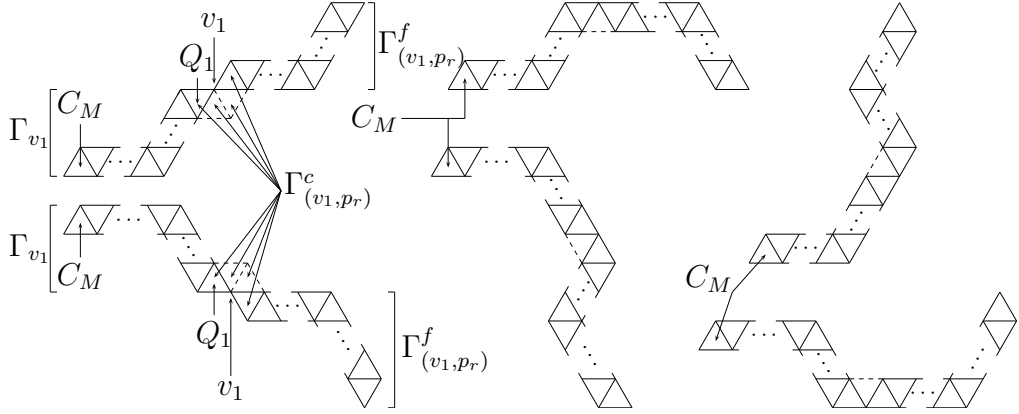}} \caption{General shapes of the
$\bar{\Gamma}_{(v_1,p_r)}$ for $SL_3$} \label{GeneralShapes}
\end{figure}

Let $\Gamma_{v_1}$ be a minimal gallery from $C_M$ to $v_1$.  As
mentioned in Section~\ref{GM}, $Q_1$ is the last alcove in
$\Gamma_{v_1}$.  Let $\Gamma^f_{(v_1,p_r)}$ be $z \Gamma_{v_1}$,
where $Q_1$ and $zQ_1$ have relative position $p_r$, $z \in SL_3
(L)$, and $z$ sends $A_M$ to $A_M$.  Let $\Gamma^c_{(v_1,p_r)}$ be
some fixed minimal connecting gallery from $Q_1$ to $zQ_1$. Define
$\bar{\Gamma}_{(v_1,p_r)} = \Gamma_{v_1} \cup \Gamma^c_{(v_1,p_r)}
\cup \Gamma^f_{(v_1,p_r)}$. The galleries
$\bar{\Gamma}_{(v_1,p_r)}$ have the general shapes pictured in
Figure~\ref{GeneralShapes}.  For clarity, only two of the
galleries in this figure have all their parts labelled.

Let $\Omega$ be a gallery in $A_M$ starting at $C_M$, and
containing any alcove at most once (so it is non-stuttering,
non-backtracking, and does not cross itself). Let
$\Omega_1,\Omega_2,\ldots,\Omega_n$ be the alcoves of $\Omega$ in
order (so $\Omega_1 = C_M$), and let $e_i$ be the edge between
$\Omega_i$ and $\Omega_{i+1}$.  Let $j$ be minimal such that $C_M$
and $\Omega_j$ are on opposite sides of the hyperplane $h_j$ in
$A_M$ determined by $e_j$.  We say $e_j$ is the first {\em choice
edge} in $\Omega$.  If $j$ does not exist then there are no choice
edges in $\Omega$. If $j$ does exist, then we define the {\em hard
choice at $e_j$} to be the gallery
$\Omega_1,\ldots,\Omega_j,\Omega_{j+1}$, and the {\em easy choice
at $e_j$} to be the gallery
$\Omega_1,\ldots,\Omega_j,f_{h_j}(\Omega_{j+1})=\Omega_j$, where
$f_{h_j}$ represents the flip of $A_M$ about $h_j$.  Given the
hard choice, we consider
$\Omega_1,\ldots,\Omega_j,\Omega_{j+1},\ldots,\Omega_n$, and we
find the minimal $k>j$ such that $h_k$ has $\Omega_k$ and $C_M$ on
opposite sides.  This is the next choice edge, given the hard
choice at $j$, and we can make easy and hard choices here.  Given
the easy choice at $j$ we consider
$\Omega_1,\ldots,\Omega_j,f_{h_j}(\Omega_{j+1}),\ldots,f_{h_j}(\Omega_{n})$,
and we find the minimal $k$ such that $k>j$, and such that
$f_{h_j}(\Omega_k)$ and $C_M$ are on opposite sides of the
hyperplane between $f_{h_j}(\Omega_k)$ and
$f_{h_j}(\Omega_{k+1})$. This is the next choice edge, given the
easy choice at $j$, and we can make either a hard or an easy
choice here. In this way we construct a binary tree, $T$, called
the {\em choice tree} for $\Omega$. Each node in $T$ except the
leaves corresponds to a choice edge in $\Omega$. At every node
except the leaves, $T$ has a branch corresponding to a hard choice
and another corresponding to an easy choice.

One can show that any non-backtracking path from the root node to
a leaf of $T$ corresponds to the retraction (onto $A_M$ centered
at $C_M$) of some gallery (or galleries) starting at $C_M$ and of
the same type at $\Omega$.  Such a path is equivalent to the
choice of a leaf of $T$, since $T$ is a tree.  The gallery
$\Omega$ itself corresponds to the path obtained by making all
hard choices in $T$. Further, all galleries starting at $C_M$ of
the same type as $\Omega$ retract in a way specified by some
non-backtracking path from the root node of $T$ to a leaf. We
define the {\em set of comprehensive folding results of $\Omega$}
to be the set of final alcoves of retractions of galleries
starting at $C_M$ that have the same type as $\Omega$. By {\em
retraction}, we always mean the retraction centered at $C_M$ onto
$A_M$. So a comprehensive folding result of $\Omega$ can also be
thought of as a non-backtracking path $F$ from the root node to a
leaf of the choice tree of $\Omega$.

We now observe that the set of comprehensive folding results of
$\Omega = \bar{\Gamma}_{(v_1,p_r)}$ contains the set of possible
$\inv(D,\sigma(D))$ for $D$ in the $(v_1,v_2,p_r)$-superpiece. We
claim that this is also an equality of sets.  Let $F$ be a
non-backtracking path from the root of $T$ to a leaf.  We define
the cf-dimension of $F$ to be
$l(\Gamma_{v_1})+l(\Gamma_{(v_1,p_r)}^c)-n_F-2$, where $n_F$ is
the number of hard choices in $F$, and $l$ represents the length
of a gallery (the number of alcoves in it). We also claim that the
cf-dimension of $F$ is equal to the dimension of the
$(v_1,v_2,p_r)$-piece of $X_{\tilde{w}}(\sigma)$, where
$\tilde{w}C_M$ is the comprehensive folding result of $\Omega$
corresponding to $F$.

We first note that one can make at most one easy choice for each
$\Omega = \bar{\Gamma}_{(v_1,p_r)}$ in Figure~\ref{GeneralShapes}.
Once this choice is made, there are no subsequent choice edges.
This can be seen just by analyzing the pictures in
Figure~\ref{GeneralShapes} on a case by case basis.  Also, one can
see that $\Gamma_{v_1} \cup \Gamma^c_{(v_1,p_r)}$ is minimal. So
the first choice edge in $\Omega$ occurs between two of the
alcoves of $\Gamma^f_{(v_1,p_r)}$.  Using these facts, one can
show that choice edges in $\Omega$ correspond to hyperplanes in
$A_M$ that pass between two alcoves of $\Gamma^f_{(v_1,p_r)}$ and
that also pass between two alcoves of $\Gamma_{v_1}$.  We now seek
to prove the claims of the previous paragraph.  Given a
non-backtracking path $F$ in $T$ from the root node to a leaf, we
need to produce some gallery $\Lambda$ such that $y\Gamma_{v_1} =
\Lambda$ for some $y \in SL_3(L)$ with $yQ_1 = Q_2$ for some $Q_2
\in \mathcal{Q}_2$, and such that $\rho_{C_M}(y^{-1}(\Lambda \cup
\Lambda^c \cup \sigma(\Lambda)))$ gives the comprehensive folding
result determined by $F$.  Here $\Lambda^c$ is a minimal gallery
from $Q_2$ to $\sigma(Q_2)$ that has the same type as
$\Gamma^c_{(v_1,p_r)}$, and $\rho_{C_M}$ is the retraction onto
$A_M$ centered at $C_M$.

Note first that $F$ determines the relative position of any two
alcoves in $\bar{\Lambda} = \Lambda \cup \Lambda^c \cup
\sigma(\Lambda)$. In our $SL_3$ case, $F$ is just an indication of
the choice edge at which to make the easy choice, if any (since
there is at most one easy choice).  We will construct $\Lambda$
starting from $\Lambda_n$, the alcove that contains $v_2$.  We
choose $\Lambda_n = Q_2$. The dimension of the set of choices for
this construction is $l(\Gamma^c_{(v_1,p_r)})-1$, since the
structure of (non-affine) Deligne-Lusztig varieties is known
\cite{DL1}.  We assume by induction that we have constructed
$\Lambda_i, \Lambda_{i+1}, \ldots , \Lambda_{n}$ (and therefore
also
$\sigma(\Lambda_n),\sigma(\Lambda_{n+1}),\ldots,\sigma(\Lambda_i)$),
and that the dimension of the space of possible such constructions
is $l(\Gamma^c_{(v_1,p_r)})+(n-i)-1-n_{(F,i)}$, where $n_{(F,i)}$
is defined as follows.  Each choice edge $e$ in $\Omega$ has two
corresponding integers $1 \leq \beta_1,\beta_2 \leq n-1$ such that
the hyperplane $h_e$ corresponding to $e$ passes between the
$\beta_1^{th}$ and $(\beta_1 +1)^{th}$ alcoves of $\Gamma_{v_1}$
(where the first alcove of $\Gamma_{v_1}$ is considered to be
$C_M$), and between the $\beta_2^{th}$ and $(\beta_2 + 1)^{th}$
alcoves of $\Gamma^f_{(v_1,p_r)}$ (where the $n^{th}$ alcove of
$\Gamma^f_{(v_1,p_r)}$ is considered to be the one containing
$v_1$). We define $n_{(F,i)}$ to be the number of choice edges $e$
such that $i \leq \beta_1,\beta_2$, and such that $F$ indicates a
hard choice at $e$. Note that if $i=1$,
$l(\Gamma^c_{(v_1,p_r)})+(n-i)-1-n_{(F,i)} = l(\Gamma_{v_1}) +
l(\Gamma^c_{(v_1,p_r)}) - n_F - 2$, and if $i=n$, then
$l(\Gamma^c_{(v_1,p_r)})+(n-i)-1-n_{(F,i)} =
l(\Gamma^c_{(v_i,p_r)})-1$.

Let $A$ be some apartment containing $\Lambda_i$ and
$\sigma(\Lambda_i)$. Let $S \subset A$ be the intersection of all
apartments that contain $\Lambda_i$ and $\sigma(\Lambda_i)$.  Let
$d_{i-1}$ be the edge of $\Lambda_i$ to which $\Lambda_{i-1}$ must
be attached (this is specified by the requirement that
$\bar{\Lambda}$  and $\Omega$ be of the same type).  Let
$\widetilde{\Lambda_{i-1}}$ be the alcove in $A$ that one gets by
reflecting $\Lambda_i$ about $d_{i-1}$.  Let
$\widetilde{\sigma(\Lambda_{i-1})}$ be the alcove in $A$ one gets
by reflecting $\sigma(\Lambda_i)$ about $\sigma(d_{i-1})$.  One
can see by considering each of the cases pictured in
Figure~\ref{GeneralShapes} that either exactly one of
$\widetilde{\Lambda_{i-1}}$ and
$\widetilde{\sigma(\Lambda_{i-1})}$ is in $S$, or neither is in
$S$.  Note that the former occurs if and only if $i-1 =
\min(\beta_1,\beta_2)$ for $\beta_1,\beta_2$ the two integers
corresponding to some choice edge in $F$.

Let $S_{i-1}$ be the intersection of all apartments containing $S
\cup \widetilde{\Lambda_{i-1}}$, and let $S^{\sigma}_{i-1}$ be the
intersection of all apartments containing $S \cup
\widetilde{\sigma(\Lambda_{i-1})}$.  One can see by considering
the cases in Figure~\ref{GeneralShapes} that if neither
$\widetilde{\Lambda_{i-1}}$ nor
$\widetilde{\sigma(\Lambda_{i-1})}$ is in $S$, then
$\widetilde{\Lambda_{i-1}}$ is not in $S_{i-1}^{\sigma}$ and
$\widetilde{\sigma(\Lambda_{i-1})}$ is not in $S_{i-1}$.
Therefore, in this case we can choose any $\Lambda_{i-1}$ adjacent
to $\Lambda_i$ by $d_{i-1}$.  This in turn determines
$\sigma(\Lambda_{i-1})$ adjacent to $\sigma(\Lambda_i)$ by
$\sigma(d_{i-1})$.  There is one dimension worth of these choices,
so the dimension of the construction down to $i-1$ is
$l(\Gamma^c_{(v_1,p_r)}) + (n-i) - 1 - n_{(F,i)} + 1$.  In this
case $n_{(F,i-1)} = n_{(F,i)}$, so
$l(\Gamma^c_{(v_1,p_r)})+(n-i)-1-n_{(F,i)}+1 =
l(\Gamma^c_{(v_1,p_r)}) + (n-(i-1))-1-n_{(F,i-1)}$.

We now consider the case in which exactly one of
$\widetilde{\Lambda_{i-1}}$ and
$\widetilde{\sigma(\Lambda_{i-1})}$ is in $S$.  We assume
$\widetilde{\Lambda_{i-1}}$ is in $S$.  The other case is similar.
This means $i-1 = \min(\beta_1,\beta_2)$ for $\beta_1,\beta_2$ the
two integers corresponding to some choice edge $e$.  If $F$
dictates a hard choice at this point, we must choose
$\Lambda_{i-1} \subset S$. There is only one such choice, causing
no increase in the dimension of the construction.  If $F$ dictates
an easy choice, we may choose any $\Lambda_{i-1}$ not in $A$, but
attached to $\Lambda_i$ via $d_{i-1}$.  There is one dimension
worth of such choices, increasing dimension by one.  In the former
case, $n_{(F,i-1)} = n_{(F,i)} + 1$, and in the latter case
$n_{(F,i+1)} = n_{(F,i)}$. In both cases, the dimension of the new
structure is $l(\Gamma^c_{(v_1,p_r)}) + (n-(i-1))-1-n_{(F,i-1)}$.
This finishes the proof of the previous claims.

The result of all this is that we can calculate the values of
$\inv(D,\sigma(D))$ for $D$ in the $(v_1,v_2,p_r)$-superpiece, and
for each $\tilde{w}$ in this set we can calculate the dimension of
the $(v_1,v_2,p_r)$-piece of $X_{\tilde{w}}(\sigma)$.  This can
all be done through straightforward computations of comprehensive
folding results and cf-dimensions.  For instance, using $v_1$ and
$p_r$ leading to the $\bar{\Gamma}_{(v_1,p_r)}$ pictured in
Figure~\ref{examplepart1}, we get the results pictured in
Figure~\ref{examplepart2}.  The numbers in
Figure~\ref{examplepart2} are the dimensions of the
$(v_1,v_2,p_r)$-pieces of the $X_{\tilde{w}}(\sigma)$, with
$\tilde{w}$ corresponding to the alcoves on which the numbers are
written. Alcoves with no numbers have empty
$(v_1,v_2,p_r)$-pieces. We did an analogous computation for every
$v_1$ in the region shown in Figure~\ref{v1regionSL3}, and for
every $p_r$ for which the corresponding $w \in W$ has $l(w) \geq
2$.  We rotated all results about the center of $C_M$ by
$120^{\circ}$ and $240^{\circ}$, combining these with the
un-rotated results.  For any alcove which contained more than one
number at that point, we took the maximum (although in all cases
for which two numbers occurred in the same alcove, these numbers
turned out to be equal). The outcome of this process was
Figure~\ref{SL3result}, which shows the $X_{\tilde{w}}(\sigma)$
that are non-empty (those corresponding to alcoves that have
numbers in them), and the dimension of these non-empty
$X_{\tilde{w}}(\sigma)$.  The bold lines in that figure correspond
to the shrunken Weyl chambers that will be discussed in
Section~\ref{formulasection}.
\begin{figure}
\centerline{\psset{xunit=0.233333333333333 in}
\psset{yunit=0.404145188432738 in}

\begin{pspicture*}(0.25,-0.25)(8.25,1.75)

  \multido{\rxb=-1.25+1,\rxt=0.75+1}{10}{ \psline[linewidth=0.05pt](\rxb,-0.25)(\rxt,1.75) }
  \multido{\rxb=0.25+1,\rxt=-1.75+1}{10}{ \psline[linewidth=0.05pt](\rxb,-0.25)(\rxt,1.75) }
  \multido{\ny=0+0.5}{4}{ \psline[linewidth=0.05pt](0,\ny)(8.25,\ny) }

\pspolygon*[linecolor=lightgray](0.5,0.5)(1.5,0.5)(1,1)
\pspolygon*[linecolor=lightgray](1.5,0.5)(1,1)(2,1)
\pspolygon*[linecolor=lightgray](2,1)(1.5,0.5)(2.5,0.5)
\pspolygon*[linecolor=lightgray](2.5,0.5)(2,1)(3,1)
\pspolygon*[linecolor=lightgray](3,1)(3.5,0.5)(4,1)
\pspolygon*[linecolor=lightgray](4,1)(3,1)(3.5,1.5)
\pspolygon*[linecolor=lightgray](2.5,0.5)(3.5,0.5)(3,1)
\pspolygon*[linecolor=lightgray](3.5,1.5)(4,1)(4.5,1.5)
\pspolygon*[linecolor=lightgray](4.5,1.5)(5,1)(5.5,1.5)
\pspolygon*[linecolor=lightgray](5,1)(5.5,1.5)(6,1)
\pspolygon*[linecolor=lightgray](5.5,1.5)(6,1)(6.5,1.5)
\pspolygon*[linecolor=lightgray](6,1)(6.5,1.5)(7,1)
\pspolygon*[linecolor=lightgray](6,1)(7,1)(6.5,0.5)
\pspolygon*[linecolor=lightgray](6.5,0.5)(7.5,0.5)(7,1)
\pspolygon*[linecolor=lightgray](7.5,0.5)(7,0)(8,0)
\pspolygon*[linecolor=lightgray](6.5,0.5)(7.5,0.5)(7,0)

\rput[bm](1,0.5){\rule[-2pt]{0pt}{5pt}$C_M$}

\end{pspicture*}}
\caption{An example of
$\bar{\Gamma}_{(v_1,p_r)}$} \label{examplepart1}
\end{figure}
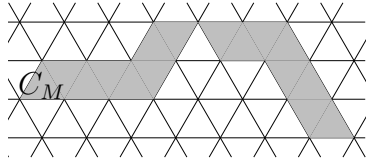
\begin{figure}
\centerline{\psset{xunit=0.233333333333333 in} \psset{yunit=0.404145188432738
in}

\begin{pspicture*}(0.25,-0.25)(8.25,2.25)

  \multido{\rxb=-2.25+1,\rxt=0.25+1}{11}{ \psline[linewidth=0.05pt](\rxb,-0.25)(\rxt,2.25) }
  \multido{\rxb=0.25+1,\rxt=-2.25+1}{11}{ \psline[linewidth=0.05pt](\rxb,-0.25)(\rxt,2.25) }
  \multido{\ny=0+0.5}{5}{ \psline[linewidth=0.05pt](0,\ny)(8.25,\ny) }

\rput[bm](1.1,0.5){\rule[-2pt]{0pt}{5pt}$C_M$}
\rput[bm](7.5,0){\rule[-2pt]{0pt}{5pt}$7$}
\rput[bm](7.5,0.65){\rule[-2pt]{0pt}{5pt}$8$}
\rput[bm](7.5,1.65){\rule[-2pt]{0pt}{5pt}$9$}

\end{pspicture*}}
\caption{Comprehensive folding
results and cf-dimensions from the example in
Figure~\ref{examplepart1}} \label{examplepart2}
\end{figure}
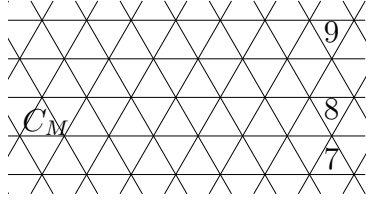
\begin{figure}
\centerline{\input{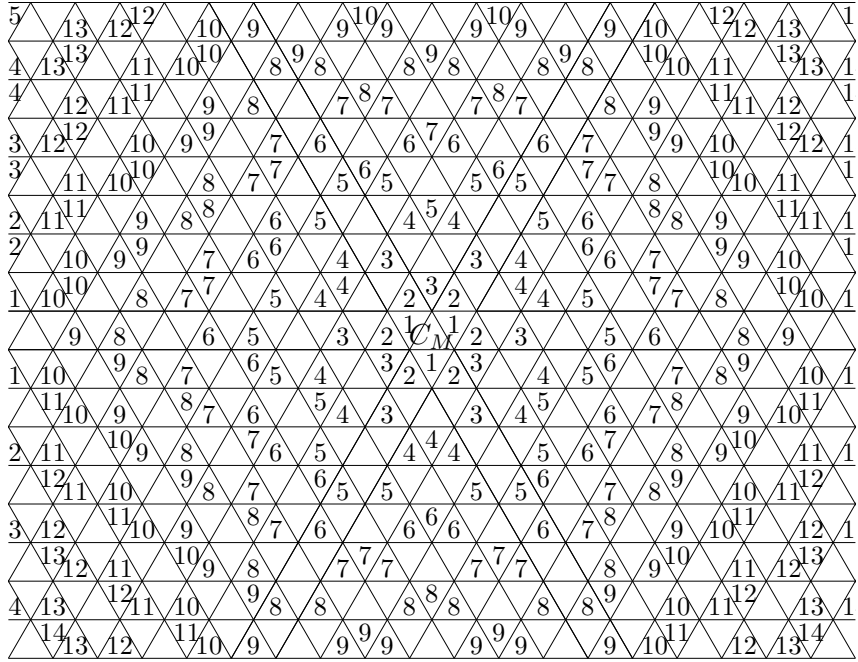}}
\caption{Main result in diagram
form for $SL_3$} \label{SL3result}
\end{figure}

Something observed in the course of the computation is that it
never happened that two different numbers occurred in the same
alcove. This means that for any fixed $\tilde{w}$, the non-empty
pieces of $X_{\tilde{w}}(\sigma)$ all have the same dimension. As
we will see in the next section, this may be related to the fact
that all vertices in the building for $SL_3$ are special.

\section{$Sp_4$}\label{Sp4}

\begin{figure}
\centerline{\psset{xunit=0.35 in} \psset{yunit=0.35 in}

\begin{pspicture*}(-0.5,-0.7)(3.7, 3.7)

    \multido{\nx=0.5+1.0}{4}{ \psline[linewidth=0.05pt](\nx,-0.7)(\nx,3.7) }
    \multido{\ny=-0.5+1.0}{5}{ \psline[linewidth=0.05pt](-0.5,\ny)(3.7,\ny) }
    \multido{\nl=-4.5+1.0,\nr=-0.3+1.0}{9}{ \psline[linewidth=0.05pt](-0.5,\nl)(3.7,\nr) }
    \multido{\nr=-4.7+1.0,\nl=-0.5+1.0}{9}{ \psline[linewidth=0.05pt](-0.5,\nl)(3.7,\nr) }

    \rput[bm](0,1.5){\rule[-2pt]{0pt}{5pt}\small{$C_M$}}
    \psset{linestyle=dashed}
    \psline[linewidth=0.05pt](-0.5,-0.7)(-0.5,3.7)

\end{pspicture*}} \caption{The region containing all
vertices $v_1$ that must be considered for $Sp_4$}
\label{v1regionSp4}
\end{figure}
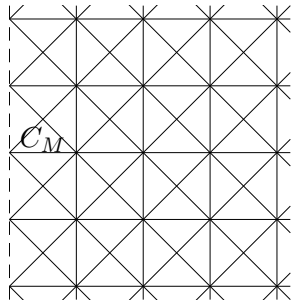

For $Sp_4$, it suffices to consider $v_1$ in the region pictured
in Figure~\ref{v1regionSp4}.  All other $v_1$ can be obtained from
these by reflecting about the line of symmetry of $C_M$.  Once
results are obtained for $v_1$ in the region specified, we will
have to reflect the results across the line of symmetry of $C_M$
as well.  Note that $v_1$ can be special or non-special for
$Sp_4$, whereas only the special case was possible for $SL_3$.

Given the restriction that $p_r$ be such that $Q_2 \cap
\mathcal{B}_1 = \{ v_2 \}$, $Q_2$ and $\sigma(Q_2)$ must share
exactly one vertex, so $p_r$ corresponds to some element of $W$ of
length $2,3$, or $4$ for $v_1$ special, and some element of length
$2$ for $v_1$ non-special.

We define $\Gamma_{v_1}$, $\Gamma^f_{(v_1,p_r)}$,
$\Gamma^c_{(v_1,p_r)}$ in the same way as in Section~\ref{SL3}.
The galleries $\bar{\Gamma}_{(v_1,p_r)}$ have the general shapes
pictured in Figure~\ref{GeneralShapesSp4nonspecial} for the case
in which $v_1$ is non-special.  For clarity, only two of the
galleries appearing in that figure have all of their parts
labelled.  Figure~\ref{GeneralShapesSp4special} contains general
shapes of the $\Gamma_{v_1}$ for $v_1$ special. The $20$ different
general shapes of the $\bar{\Gamma}_{(v_1,p_r)}$ can be deduced
from these $4$ possible $\Gamma_{v_1}$ by determining
$\Gamma^f_{(v_1,p_r)}$ and $\Gamma^c_{(v_1,p_r)}$ from each
$\Gamma_{v_1}$ using each of the $5$ possible $p_r$.
\begin{figure}
\centerline{\psset{xunit=0.3 in} \psset{yunit=0.3 in}

\begin{pspicture*}(-0.5,-0.6)(18.5,9)


    \psline[linewidth=0.05pt](-0.5,-0.5)(0.5,0.5)
    \psline[linewidth=0.05pt](0.5,0.5)(1,0)
    \psline[linewidth=0.05pt](1,0)(0.5,-0.5)
    \psline[linewidth=0.05pt](0.5,-0.5)(-0.5,-0.5)
    \psline[linewidth=0.05pt](0.5,-0.5)(0,0)
    \psline[linewidth=0.05pt](0.5,-0.5)(0.5,0.5)
    \psline[linewidth=0.05pt](1.5,0.5)(1.5,2.5)
    \psline[linewidth=0.05pt](1.5,2.5)(0.5,2.5)
    \psline[linewidth=0.05pt](0.5,2.5)(1.5,1.5)
    \psline[linewidth=0.05pt](1.5,2.5)(0.5,1.5)
    \psline[linewidth=0.05pt](0.5,1.5)(1.5,1.5)
    \psline[linewidth=0.05pt](0.5,1.5)(1.5,0.5)
    \psline[linewidth=0.05pt](1,1)(1.5,1.5)
    \psline[linewidth=0.05pt](0.5,3.5)(1.5,3.5)
    \psline[linewidth=0.05pt](1.5,3.5)(0.5,4.5)
    \psline[linewidth=0.05pt](0.5,4.5)(1.5,4.5)
    \psline[linewidth=0.05pt](1.5,4.5)(0.5,3.5)
    \psline[linewidth=0.05pt](0.5,5.5)(0.5,7.5)
    \psline[linewidth=0.05pt](0.5,7.5)(1.5,6.5)
    \psline[linewidth=0.05pt](1.5,6.5)(0.5,5.5)
    \psline[linewidth=0.05pt](0.5,5.5)(1.5,5.5)
    \psline[linewidth=0.05pt](1.5,5.5)(0.5,6.5)
    \psline[linewidth=0.05pt](0.5,6.5)(1.5,6.5)
    \psline[linewidth=0.05pt](0.5,6.5)(1,7)
    \psline[linewidth=0.05pt](1.5,7.5)(1,8)
    \psline[linewidth=0.05pt](1,8)(1.5,8.5)
    \psline[linewidth=0.05pt](1.5,7.5)(1.5,8.5)
    \psline[linewidth=0.05pt](1.5,8.5)(2.5,8.5)
    \psline[linewidth=0.05pt](2.5,8.5)(1.5,7.5)
    \psline[linewidth=0.05pt](2,8)(1.5,8.5)
    \psline[linewidth=0.05pt](2.8,4.05)(2.8,8.5)
    \psline[linewidth=0.05pt](2.6,4.05)(2.8,4.05)
    \psline[linewidth=0.05pt](2.6,8.5)(2.8,8.5)
    \rput[ml](2.9,1){$\Gamma_{v_1}$}
    \psline[linewidth=0.05pt](2.8,-0.5)(2.8,3.95)
    \psline[linewidth=0.05pt](2.6,-0.5)(2.8,-0.5)
    \psline[linewidth=0.05pt](2.6,3.95)(2.8,3.95)
    \rput[ml](2.9,6.25){$\Gamma^f_{(v_1,p_r)}$}
    \psset{linestyle=dashed}
        \psline[linewidth=0.05pt](1.5,3.5)(1.5,4.5)
    \psset{linestyle=solid}
    \psline[linewidth=0.05pt]{->}(3.5,2.5)(1,3.667)
    \psline[linewidth=0.05pt]{->}(3.5,2.5)(1,4.333)
    \psline[linewidth=0.05pt]{->}(3.5,2.5)(1.333,3.9)
    \rput[ml](3.6,2.5){$\Gamma^c_{(v_1,p_r)}$}
    \rput[c](0,-0.333){$C_M$}
    \psline[linewidth=0.05pt]{->}(0.2,3.667)(0.9,3.667)
    \rput[mr](0.1,3.667){$Q_1$}
    \rput[mr](0.1,4.333){$v_1$}
    \psline[linewidth=0.05pt]{->}(0.1,4.333)(1,4)
    \rput[c]{45}(1,0.5){\scriptsize{$\cdots$}}
    \rput[c]{90}(1,3){\scriptsize{$\cdots$}}
    \rput[c]{90}(1,5){\scriptsize{$\cdots$}}
    \rput[c]{45}(1,7.5){\scriptsize{$\cdots$}}

    \psline[linewidth=0.05pt](7.5,0.5)(8.5,0.5)
    \psline[linewidth=0.05pt](7.5,0.5)(8,1)
    \psline[linewidth=0.05pt](8,1)(8.5,0.5)
    \psline[linewidth=0.05pt](8.5,0.5)(8,0)
    \psline[linewidth=0.05pt](8,0)(7.5,0.5)
    \psline[linewidth=0.05pt](7,1)(6.5,1.5)
    \psline[linewidth=0.05pt](6.5,1.5)(6.5,2.5)
    \psline[linewidth=0.05pt](6.5,2.5)(7.5,1.5)
    \psline[linewidth=0.05pt](6.5,2.5)(7.5,2.5)
    \psline[linewidth=0.05pt](7.5,2.5)(6.5,1.5)
    \psline[linewidth=0.05pt](6.5,1.5)(7.5,1.5)
    \psline[linewidth=0.05pt](7.5,1.5)(7,1)
    \psline[linewidth=0.05pt](6.5,3.5)(7.5,4.5)
    \psline[linewidth=0.05pt](7.5,4.5)(6.5,4.5)
    \psline[linewidth=0.05pt](6.5,4.5)(7.5,3.5)
    \psline[linewidth=0.05pt](7.5,3.5)(6.5,3.5)
    \psline[linewidth=0.05pt](6.5,5.5)(7.5,6.5)
    \psline[linewidth=0.05pt](7.5,6.5)(6.5,6.5)
    \psline[linewidth=0.05pt](6.5,6.5)(7.5,5.5)
    \psline[linewidth=0.05pt](7.5,5.5)(6.5,5.5)
    \psline[linewidth=0.05pt](6.5,6.5)(7,7)
    \psline[linewidth=0.05pt](7,7)(7.5,6.5)
    \psline[linewidth=0.05pt](7.5,6.5)(7.5,5.5)
    \psline[linewidth=0.05pt](5.5,7.5)(6.5,7.5)
    \psline[linewidth=0.05pt](6.5,7.5)(6,8)
    \psline[linewidth=0.05pt](6,8)(5.5,7.5)
    \psline[linewidth=0.05pt](5.5,7.5)(6,7)
    \psline[linewidth=0.05pt](6,7)(6.5,7.5)
    \psline[linewidth=0.05pt](5.4,8)(5.4,4.05)
    \psline[linewidth=0.05pt](5.6,4.05)(5.4,4.05)
    \psline[linewidth=0.05pt](5.6,8)(5.4,8)
    \rput[mr](5.4,5){$\Gamma_{v_1}$}
    \psline[linewidth=0.05pt](5.4,3.95)(5.4,0)
    \psline[linewidth=0.05pt](5.4,3.95)(5.6,3.95)
    \psline[linewidth=0.05pt](5.4,0)(5.6,0)
    \rput[mr](5.4,1){$\Gamma^f_{(v_1,p_r)}$}
    \psset{linestyle=dashed}
        \psline[linewidth=0.05pt](6.5,3.5)(6.5,4.5)
    \psset{linestyle=solid}
    \psline[linewidth=0.05pt]{->}(4.9,2.5)(7,3.667)
    \psline[linewidth=0.05pt]{->}(4.9,2.5)(7,4.333)
    \psline[linewidth=0.05pt]{->}(4.9,2.5)(6.667,3.9)
    \rput[c](6,7.667){$C_M$}
    \psline[linewidth=0.05pt]{->}(7.8,3.667)(7.1,3.667)
    \rput[ml](7.9,3.667){$Q_1$}
    \rput[ml](7.9,4.333){$v_1$}
    \psline[linewidth=0.05pt]{->}(7.8,4.333)(7,4)
    \rput[c]{-45}(6.5,7){\scriptsize{$\cdots$}}
    \rput[c]{90}(7,5){\scriptsize{$\cdots$}}
    \rput[c]{90}(7,3){\scriptsize{$\cdots$}}
    \rput[c]{-45}(7.5,1){\scriptsize{$\cdots$}}

    \psline[linewidth=0.05pt](9.5,0.5)(10.5,1.5)
    \psline[linewidth=0.05pt](10.5,1.5)(11,1)
    \psline[linewidth=0.05pt](11,1)(10.5,0.5)
    \psline[linewidth=0.05pt](10.5,1.5)(10.5,0.5)
    \psline[linewidth=0.05pt](10.5,0.5)(10,1)
    \psline[linewidth=0.05pt](10.5,0.5)(9.5,0.5)
    \psline[linewidth=0.05pt](11,2)(11.5,1.5)
    \psline[linewidth=0.05pt](11,2)(11.5,2.5)
    \psline[linewidth=0.05pt](11.5,1.5)(11.5,2.5)
    \psline[linewidth=0.05pt](11.5,2.5)(12.5,1.5)
    \psline[linewidth=0.05pt](12.5,1.5)(12.5,2.5)
    \psline[linewidth=0.05pt](12.5,2.5)(11.5,1.5)
    \psline[linewidth=0.05pt](13.5,2.5)(13.5,1.5)
    \psline[linewidth=0.05pt](13.5,1.5)(14.5,2.5)
    \psline[linewidth=0.05pt](14.5,2.5)(14.5,1.5)
    \psline[linewidth=0.05pt](14.5,1.5)(13.5,2.5)
    \psline[linewidth=0.05pt](11.5,1.5)(12.5,1.5)
    \psline[linewidth=0.05pt](15.5,1.5)(15.5,2.5)
    \psline[linewidth=0.05pt](15.5,2.5)(16.5,1.5)
    \psline[linewidth=0.05pt](16.5,1.5)(16.5,2.5)
    \psline[linewidth=0.05pt](16.5,2.5)(15.5,1.5)
    \psline[linewidth=0.05pt](15.5,2.5)(16.5,2.5)
    \psline[linewidth=0.05pt](16.5,2.5)(17,2)
    \psline[linewidth=0.05pt](17,2)(16.5,1.5)
    \psline[linewidth=0.05pt](17,3)(17.5,2.5)
    \psline[linewidth=0.05pt](17.5,2.5)(17.5,3.5)
    \psline[linewidth=0.05pt](17.5,3.5)(17,3)
    \psline[linewidth=0.05pt](17.5,2.5)(18.5,3.5)
    \psline[linewidth=0.05pt](17.5,3.5)(18.5,3.5)
    \psline[linewidth=0.05pt](17.5,3.5)(18,3)
    \rput[c](10,0.667){$C_M$}
    \rput[c]{45}(11,1.5){\scriptsize{$\cdots$}}
    \rput[c](13,2){\scriptsize{$\cdots$}}
    \rput[c](15,2){\scriptsize{$\cdots$}}
    \rput[c]{45}(17,2.5){\scriptsize{$\cdots$}}
    \psset{linestyle=dashed}
        \psline[linewidth=0.05pt](13.5,2.5)(14.5,2.5)
    \psset{linestyle=solid}

    \psline[linewidth=0.05pt](8.5,7.5)(9,8)
    \psline[linewidth=0.05pt](9,8)(9.5,7.5)
    \psline[linewidth=0.05pt](9.5,7.5)(9,7)
    \psline[linewidth=0.05pt](9,7)(8.5,7.5)
    \psline[linewidth=0.05pt](8.5,7.5)(9.5,7.5)
    \psline[linewidth=0.05pt](9.5,6.5)(10.5,7.5)
    \psline[linewidth=0.05pt](10.5,7.5)(11.5,6.5)
    \psline[linewidth=0.05pt](11.5,6.5)(9.5,6.5)
    \psline[linewidth=0.05pt](11.5,6.5)(11.5,7.5)
    \psline[linewidth=0.05pt](11.5,7.5)(10.5,6.5)
    \psline[linewidth=0.05pt](10.5,6.5)(10.5,7.5)
    \psline[linewidth=0.05pt](10.5,6.5)(10,7)
    \psline[linewidth=0.05pt](12.5,6.5)(12.5,7.5)
    \psline[linewidth=0.05pt](12.5,7.5)(13.5,6.5)
    \psline[linewidth=0.05pt](13.5,6.5)(13.5,7.5)
    \psline[linewidth=0.05pt](13.5,7.5)(12.5,6.5)
    \psline[linewidth=0.05pt](14.5,6.5)(14.5,7.5)
    \psline[linewidth=0.05pt](14.5,7.5)(15.5,6.5)
    \psline[linewidth=0.05pt](15.5,6.5)(16.5,7.5)
    \psline[linewidth=0.05pt](16.5,7.5)(14.5,7.5)
    \psline[linewidth=0.05pt](14.5,6.5)(15.5,7.5)
    \psline[linewidth=0.05pt](15.5,7.5)(15.5,6.5)
    \psline[linewidth=0.05pt](15.5,7.5)(16,7)
    \psline[linewidth=0.05pt](16.5,6.5)(17,7)
    \psline[linewidth=0.05pt](17,7)(17.5,6.5)
    \psline[linewidth=0.05pt](17.5,6.5)(16.5,6.5)
    \psline[linewidth=0.05pt](17.5,6.5)(17,6)
    \psline[linewidth=0.05pt](17,6)(16.5,6.5)
    \rput[c](9,7.667){$C_M$}
    \rput[c]{-45}(9.5,7){\scriptsize{$\cdots$}}
    \rput[c](12,7){\scriptsize{$\cdots$}}
    \rput[c](14,7){\scriptsize{$\cdots$}}
    \rput[c]{-45}(16.5,7){\scriptsize{$\cdots$}}
    \psset{linestyle=dashed}
        \psline[linewidth=0.05pt](12.5,6.5)(13.5,6.5)
    \psset{linestyle=solid}

\end{pspicture*}}
\caption{General shapes of the
$\bar{\Gamma}_{(v_1,p_r)}$ for $Sp_4$, non-special $v_1$}
\label{GeneralShapesSp4nonspecial}
\end{figure}
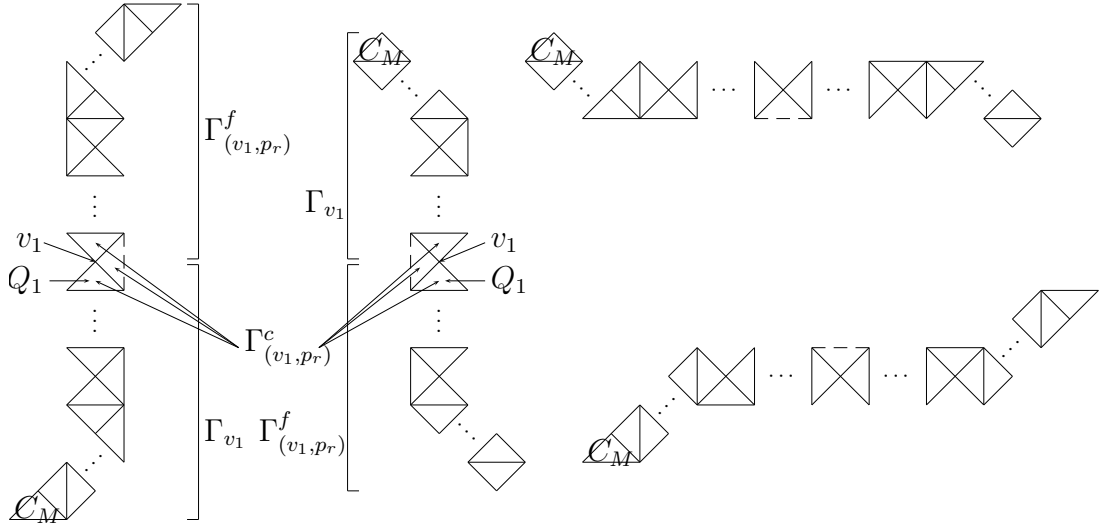
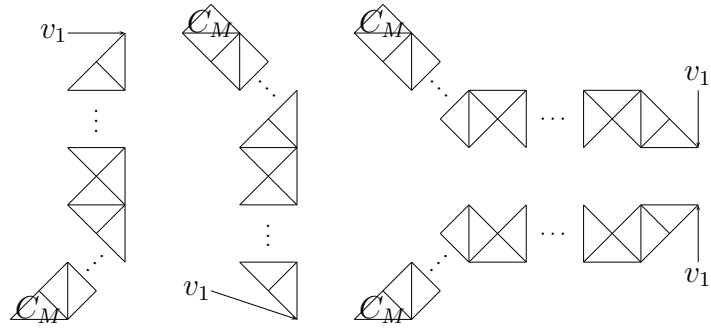
\begin{figure}
\centerline{\psset{xunit=0.3 in} \psset{yunit=0.3 in}

\begin{pspicture*}(-0.5,-0.7)(11.8,5)

    \psline[linewidth=0.05pt](-0.5,-0.5)(0.5,0.5)
    \psline[linewidth=0.05pt](0.5,0.5)(0.5,-0.5)
    \psline[linewidth=0.05pt](0.5,-0.5)(0,0)
    \psline[linewidth=0.05pt](0.5,-0.5)(-0.5,-0.5)
    \psline[linewidth=0.05pt](1.5,0.5)(1.5,2.5)
    \psline[linewidth=0.05pt](1.5,2.5)(0.5,1.5)
    \psline[linewidth=0.05pt](0.5,1.5)(1.5,0.5)
    \psline[linewidth=0.05pt](1.5,2.5)(0.5,2.5)
    \psline[linewidth=0.05pt](1.5,1.5)(1,1)
    \psline[linewidth=0.05pt](1.5,1.5)(0.5,1.5)
    \psline[linewidth=0.05pt](1.5,1.5)(0.5,2.5)
    \psline[linewidth=0.05pt](0.5,-0.5)(1,0)
    \psline[linewidth=0.05pt](1,0)(0.5,0.5)
    \psline[linewidth=0.05pt](0.5,3.5)(1.5,3.5)
    \psline[linewidth=0.05pt](1.5,3.5)(1,4)
    \psline[linewidth=0.05pt](1.5,3.5)(1.5,4.5)
    \psline[linewidth=0.05pt](1.5,4.5)(0.5,3.5)
    \rput[c](0,-0.333){$C_M$}
    \rput[mr](0.5,4.5){$v_1$}
    \psline[linewidth=0.05pt]{->}(0.5,4.5)(1.5,4.5)
    \rput[c]{45}(1,0.5){\scriptsize{$\cdots$}}
    \rput[c]{90}(1,3){\scriptsize{$\cdots$}}

    \psline[linewidth=0.05pt](2.5,4.5)(3,5)
    \psline[linewidth=0.05pt](3,5)(4,4)
    \psline[linewidth=0.05pt](4,4)(3.5,3.5)
    \psline[linewidth=0.05pt](3.5,3.5)(2.5,4.5)
    \psline[linewidth=0.05pt](3.5,3.5)(3.5,4.5)
    \psline[linewidth=0.05pt](3,4)(3.5,4.5)
    \psline[linewidth=0.05pt](2.5,4.5)(3.5,4.5)
    \psline[linewidth=0.05pt](3.5,2.5)(4.5,3.5)
    \psline[linewidth=0.05pt](4.5,3.5)(4.5,1.5)
    \psline[linewidth=0.05pt](4.5,1.5)(3.5,2.5)
    \psline[linewidth=0.05pt](4.5,2.5)(3.5,1.5)
    \psline[linewidth=0.05pt](3.5,1.5)(4.5,1.5)
    \psline[linewidth=0.05pt](4.5,2.5)(3.5,2.5)
    \psline[linewidth=0.05pt](4.5,2.5)(4,3)
    \psline[linewidth=0.05pt](3.5,0.5)(4.5,-0.5)
    \psline[linewidth=0.05pt](4.5,-0.5)(4.5,0.5)
    \psline[linewidth=0.05pt](4.5,0.5)(3.5,0.5)
    \psline[linewidth=0.05pt](4.5,0.5)(4,0)
    \rput[c](3,4.667){$C_M$}
    \rput[mr](3,0){$v_1$}
    \psline[linewidth=0.05pt]{->}(3,0)(4.5,-0.5)
    \rput[c]{-45}(4,3.5){\scriptsize{$\cdots$}}
    \rput[c]{90}(4,1){\scriptsize{$\cdots$}}

    \psline[linewidth=0.05pt](5.5,-0.5)(6.5,0.5)
    \psline[linewidth=0.05pt](6.5,0.5)(7,0)
    \psline[linewidth=0.05pt](7,0)(6.5,-0.5)
    \psline[linewidth=0.05pt](6.5,-0.5)(5.5,-0.5)
    \psline[linewidth=0.05pt](6.5,-0.5)(6.5,0.5)
    \psline[linewidth=0.05pt](6.5,-0.5)(6,0)
    \psline[linewidth=0.05pt](7,1)(7.5,1.5)
    \psline[linewidth=0.05pt](7,1)(7.5,0.5)
    \psline[linewidth=0.05pt](7.5,0.5)(7.5,1.5)
    \psline[linewidth=0.05pt](7.5,1.5)(8.5,0.5)
    \psline[linewidth=0.05pt](8.5,0.5)(8.5,1.5)
    \psline[linewidth=0.05pt](8.5,1.5)(7.5,0.5)
    \psline[linewidth=0.05pt](7.5,0.5)(8.5,0.5)
    \psline[linewidth=0.05pt](9.5,0.5)(9.5,1.5)
    \psline[linewidth=0.05pt](9.5,1.5)(10.5,0.5)
    \psline[linewidth=0.05pt](10.5,0.5)(10.5,1.5)
    \psline[linewidth=0.05pt](10.5,1.5)(9.5,0.5)
    \psline[linewidth=0.05pt](9.5,0.5)(10.5,0.5)
    \psline[linewidth=0.05pt](10.5,0.5)(11.5,1.5)
    \psline[linewidth=0.05pt](11.5,1.5)(10.5,1.5)
    \psline[linewidth=0.05pt](11,1)(10.5,1.5)
    \rput[c](6,-0.333){$C_M$}
    \rput[mt](11.5,0.4){$v_1$}
    \psline[linewidth=0.05pt]{->}(11.5,0.5)(11.5,1.5)
    \rput[c]{45}(7,0.5){\scriptsize{$\cdots$}}
    \rput[c](9,1){\scriptsize{$\cdots$}}

    \psline[linewidth=0.05pt](5.5,4.5)(6,5)
    \psline[linewidth=0.05pt](6,5)(7,4)
    \psline[linewidth=0.05pt](7,4)(6.5,3.5)
    \psline[linewidth=0.05pt](6.5,3.5)(5.5,4.5)
    \psline[linewidth=0.05pt](6.5,4.5)(5.5,4.5)
    \psline[linewidth=0.05pt](6.5,4.5)(6.5,3.5)
    \psline[linewidth=0.05pt](6.5,4.5)(6,4)
    \psline[linewidth=0.05pt](7,3)(7.5,3.5)
    \psline[linewidth=0.05pt](7,3)(7.5,2.5)
    \psline[linewidth=0.05pt](7.5,2.5)(7.5,3.5)
    \psline[linewidth=0.05pt](7.5,3.5)(8.5,2.5)
    \psline[linewidth=0.05pt](8.5,2.5)(8.5,3.5)
    \psline[linewidth=0.05pt](8.5,3.5)(7.5,2.5)
    \psline[linewidth=0.05pt](7.5,3.5)(8.5,3.5)
    \psline[linewidth=0.05pt](9.5,2.5)(9.5,3.5)
    \psline[linewidth=0.05pt](9.5,3.5)(10.5,2.5)
    \psline[linewidth=0.05pt](10.5,2.5)(10.5,3.5)
    \psline[linewidth=0.05pt](10.5,3.5)(9.5,2.5)
    \psline[linewidth=0.05pt](10.5,3.5)(9.5,3.5)
    \psline[linewidth=0.05pt](10.5,3.5)(11.5,2.5)
    \psline[linewidth=0.05pt](11.5,2.5)(10.5,2.5)
    \psline[linewidth=0.05pt](10.5,2.5)(11,3)
    \rput[c](6,4.667){$C_M$}
    \rput[mb](11.5,3.6){$v_1$}
    \psline[linewidth=0.05pt]{->}(11.5,3.5)(11.5,2.5)
    \rput[c]{-45}(7,3.5){\scriptsize{$\cdots$}}
    \rput[c](9,3){\scriptsize{$\cdots$}}

\end{pspicture*}}
\caption{General shapes of the
$\Gamma_{v_1}$ for $Sp_4$, special $v_1$}
\label{GeneralShapesSp4special}
\end{figure}

We define {\em choice edges}, {\em hard} and {\em easy choices},
and the {\em choice tree} of a non-stuttering, non-backtracking
gallery $\Omega$ in $A_M$ that does not cross itself just as we
did for $SL_3$. Again call the choice tree $T$. As for $SL_3$, any
non-backtracking path from the root node of $T$ to a leaf
corresponds to the retraction of some gallery starting at $C_M$ of
the same type as $\Omega$.  All galleries starting at $C_M$ of the
same type as $\Omega$ retract in a way specified by some
non-backtracking path from the root node of $T$ to a leaf. We
define the set of comprehensive folding results of $\Omega$ as
before.

The set of comprehensive folding results of $\Omega =
\bar{\Gamma}_{(v_1,p_r)}$ contains the set of possible
$\inv(D,\sigma(D))$ for $D$ in the $(v_1,v_2,p_r)$-superpiece. We
claim that this is an equality of sets.  We define cf-dimension of
a non-backtracking path $F$ from root to leaf as before.  We claim
that the cf-dimension of $F$ is equal to the dimension of the
$(v_1,v_2,p_r)$-piece of $X_{\tilde{w}}(\sigma)$, where
$\tilde{w}C_M$ is the comprehensive folding result of $\Omega$
corresponding to $F$.

We first note that $F$ can contain at most two easy choices. In
fact, the maximum number of easy choices that $F$ can contain is
$-m+4$, where $m$ is the length of $p_r$ in $W$.  This result is
obtained by considering cases.  For $SL_3$, the maximum number of
easy choices is $-m+3$.  As in the $SL_3$ case, for $Sp_4$,
$\Gamma_{v_1} \cup \Gamma^c_{(v_1,p_r)}$ is minimal.

We define a {\em non-primal choice edge} to be a node in $T$ that
occurs below some easy choice in $T$ (i.e., the non-backtracking
path from the root node to the node in question passes through an
edge in $T$ corresponding to an easy choice).  A {\em primal
choice edge} is any choice edge that is not non-primal. All choice
edges for $SL_3$ are primal.  For $SL_3$ and $Sp_4$, all primal
choice edges in $\Omega=\bar{\Gamma}_{(v_1,p_r)}$ correspond to
hyperplanes in $A_M$ that pass between two alcoves of
$\Gamma^f_{(v_1,p_r)}$ and that also pass between two alcoves of
$\Gamma_{v_1}$.

Given a primal choice edge in $F$, we define the two corresponding
integers $1 \leq \beta_1,\beta_2 \leq n-1$ as in the $SL_3$ case.
We will also define $1 \leq \beta_1,\beta_2 \leq n-1$ for a
non-primal choice edge $e$, but in a slightly different way. Since
$e$ is non-primal, there is some choice edge $d$ above $e$ in $F$
at which $F$ makes the easy choice. Let $h_{d}$ be the hyperplane
in $A_M$ determined by the edge $d$ in $\Omega$ (so $h_{d}$ is a
hyperplane separating two alcoves of $\Gamma^f_{(v_1,p_r)}$ and
also two alcoves of $\Gamma_{v_1}$). Let $f_{h_{d}}$ be the flip
in $A_M$ about $h_{d}$. Consider the gallery in $A_M$ obtained by
applying $f_{h_{d}}$ to the alcoves in $\Omega$ that occur after
$d$ (here $C_M$ is considered to be the first alcove of $\Omega$).
Let $\tilde{e}= f_{h_{d}}(e)$, and let $h_{\tilde{e}}$ be the
hyperplane in $A_M$ determined by $\tilde{e}$.  Let $\beta_2$ be
such that $h_{e}$ passes between the $\beta_2^{th}$ and
$(\beta_2+1)^{th}$ alcoves of $\Gamma^f_{(v_1,p_r)}$ (here the
$n^{th}$ alcove of $\Gamma^f_{(v_1,p_r)}$ is considered to be the
one containing $v_1$). If $h_{\tilde{e}}$ passes between two
alcoves of $\Gamma_{v_1}$, then let $\beta_1$ be such that
$h_{\tilde{e}}$ passes between the $\beta_1^{th}$ and
$(\beta_1+1)^{th}$ alcoves of $\Gamma_{v_1}$ (here $C_M$ is
considered to be the first alcove of $\Gamma_{v_1})$.  Otherwise
let $\beta_1 = n$.

Now, given a non-backtracking path $F$ from the root node of $T$
to a leaf, we want to produce a gallery $\Lambda$ such that
$y\Gamma_{v_1}=\Lambda$ for some $y \in Sp_4(L)$ with $yQ_1=Q_2$
for some $Q_2 \in \mathcal{Q}_2$, and such that
$\rho_{C_M}(y^{-1}(\Lambda \cup \Lambda^c \cup \sigma(\Lambda)))$
gives the comprehensive folding result determined by $F$.  Here,
as before, $\Lambda^c$ is a minimal gallery from $Q_2$ to
$\sigma(Q_2)$ that has the same type as $\Gamma^c_{(v_1,p_r)}$.

We choose $\Lambda_n = Q_2$.  The dimension of the set of such
choices is $l(\Gamma^c_{(v_1,p_r)})-1$ \cite{DL1}.  We assume by
induction that we have constructed $\Lambda_i,\ldots,\Lambda_n$
and $\sigma(\Lambda_n),\ldots,\sigma(\Lambda_i)$, and that the
dimension of the space of choices for this construction is
$l(\Gamma^c_{(v_1,p_r)})+(n-i)-1-n_{(F,i)}$, where $n_{(F,i)}$ is
defined to be the number of choice edges $e$ in $F$ such that $i
\leq \beta_1,\beta_2$ and such that $F$ indicates a hard choice at
$e$. Here $\beta_1,\beta_2$ are the integers corresponding to $e$,
defined in the new way. If $i=1$,
$l(\Gamma^c_{(v_1,p_r)})+(n-i)-1-n_{(F,i)} =
l(\Gamma_{v_1})+l(\Gamma^c_{(v_1,p_r)})-n_F-2$, and if $i=n$ then
$l(\Gamma^c_{(v_1,p_r)})+(n-i)-1-n_{(F,i)} =
l(\Gamma^c_{(v_1,p_r)})-1$.

As was the case for $SL_3$, let $A$ be some apartment containing
$\Lambda_i$ and $\sigma(\Lambda_i)$.  Let $S \subset A$ be the
intersection of all apartments that contain both $\Lambda_i$ and
$\sigma(\Lambda_i)$.  Let $d_{i-1}$ be the edge of $\Lambda_i$ to
which $\Lambda_{i-1}$ must be attached.  Let
$\widetilde{\Lambda_{i-1}}$ be the alcove in $A$ that one gets by
reflecting $\Lambda_i$ about $d_{i-1}$.  Let
$\widetilde{\sigma(\Lambda_{i-1})}$ be the alcove in $A$ one gets
by reflecting $\sigma(\Lambda_i)$ about $\sigma(d_{i-1})$.  Let
$S_{i-1}$ be the intersection of all apartments containing $S \cup
\widetilde{\Lambda_{i-1}}$ and let $S^{\sigma}_{i-1}$ be the
intersection of all apartments containing $S \cup
\widetilde{\sigma(\Lambda_{i-1})}$.

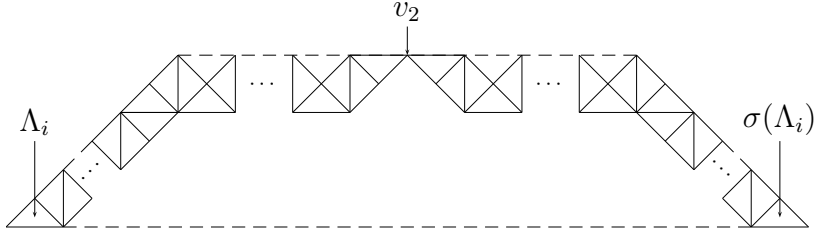
\begin{figure}
\centerline{\psset{xunit=0.3 in} \psset{yunit=0.3 in}

\begin{pspicture*}(-0.5,-0.7)(13.8,3.5)

    \psline[linewidth=0.05pt](-0.5,-0.5)(0.5,0.5)
    \psline[linewidth=0.05pt](0.5,0.5)(1,0)
    \psline[linewidth=0.05pt](1,0)(0.5,-0.5)
    \psline[linewidth=0.05pt](0.5,-0.5)(-0.5,-0.5)
    \psline[linewidth=0.05pt](0.5,-0.5)(0.5,0.5)
    \psline[linewidth=0.05pt](0,0)(0.5,-0.5)
    \psline[linewidth=0.05pt](1,1)(1.5,0.5)
    \psline[linewidth=0.05pt](1.5,0.5)(3.5,2.5)
    \psline[linewidth=0.05pt](3.5,2.5)(3.5,1.5)
    \psline[linewidth=0.05pt](3.5,1.5)(2.5,2.5)
    \psline[linewidth=0.05pt](3.5,1.5)(1.5,1.5)
    \psline[linewidth=0.05pt](2.5,2.5)(1,1)
    \psline[linewidth=0.05pt](2.5,1.5)(2,2)
    \psline[linewidth=0.05pt](2.5,1.5)(1.5,1.5)
    \psline[linewidth=0.05pt](1.5,1.5)(2,1)
    \psline[linewidth=0.05pt](1.5,1.5)(1.5,0.5)
    \psline[linewidth=0.05pt](2.5,1.5)(2.5,2.5)
    \psline[linewidth=0.05pt](4.5,2.5)(4.5,1.5)
    \psline[linewidth=0.05pt](4.5,1.5)(5.5,2.5)
    \psline[linewidth=0.05pt](5.5,2.5)(5.5,1.5)
    \psline[linewidth=0.05pt](5.5,1.5)(4.5,2.5)
    \psline[linewidth=0.05pt](5.5,1.5)(4.5,1.5)
    \psline[linewidth=0.05pt](5.5,1.5)(6.5,2.5)
    \psline[linewidth=0.05pt](6.5,2.5)(7.5,1.5)
    \psline[linewidth=0.05pt](5.5,2.5)(6,2)
    \psline[linewidth=0.05pt](7.5,1.5)(7.5,2.5)
    \psline[linewidth=0.05pt](7.5,2.5)(5.5,2.5)
    \psline[linewidth=0.05pt](7.5,2.5)(7,2)
    \psline[linewidth=0.05pt](9.5,1.5)(9.5,2.5)
    \psline[linewidth=0.05pt](9.5,2.5)(11.5,0.5)
    \psline[linewidth=0.05pt](10.5,1.5)(10.5,2.5)
    \psline[linewidth=0.05pt](10.5,2.5)(9.5,1.5)
    \psline[linewidth=0.05pt](9.5,1.5)(11.5,1.5)
    \psline[linewidth=0.05pt](11.5,0.5)(12,1)
    \psline[linewidth=0.05pt](12,1)(10.5,2.5)
    \psline[linewidth=0.05pt](11.5,1.5)(11.5,0.5)
    \psline[linewidth=0.05pt](11.5,1.5)(11,1)
    \psline[linewidth=0.05pt](10.5,1.5)(11,2)
    \psline[linewidth=0.05pt](7.5,2.5)(8.5,1.5)
    \psline[linewidth=0.05pt](8.5,1.5)(8.5,2.5)
    \psline[linewidth=0.05pt](8.5,2.5)(7.5,1.5)
    \psline[linewidth=0.05pt](7.5,1.5)(8.5,1.5)
    \psline[linewidth=0.05pt](12,0)(12.5,0.5)
    \psline[linewidth=0.05pt](12.5,0.5)(12.5,-0.5)
    \psline[linewidth=0.05pt](12.5,-0.5)(12,0)
    \psline[linewidth=0.05pt](13.5,-0.5)(12.5,0.5)
    \psline[linewidth=0.05pt](13.5,-0.5)(12.5,-0.5)
    \psline[linewidth=0.05pt](12.5,-0.5)(13,0)

    \psset{linestyle=dashed}
        \psline[linewidth=0.05pt](2.5,2.5)(10.5,2.5)
        \psline[linewidth=0.05pt](0.5,0.5)(1,1)
        \psline[linewidth=0.05pt](0.5,-0.5)(12.5,-0.5)
        \psline[linewidth=0.05pt](12,1)(12.5,0.5)
    \psset{linestyle=solid}

    \rput[mb](0,1.1){$\Lambda_i$}
    \psline[linewidth=0.05pt]{->}(0,1)(0,-0.333)
    \rput[mb](13,1.1){$\sigma(\Lambda_i)$}
    \psline[linewidth=0.05pt]{->}(13,1)(13,-0.333)
    \rput[c]{45}(1,0.5){\scriptsize{$\cdots$}}
    \rput[c](4,2){\scriptsize{$\cdots$}}
    \rput[c](9,2){\scriptsize{$\cdots$}}
    \rput[c]{-45}(12,0.5){\scriptsize{$\cdots$}}
    \rput[mb](6.5,3.1){$v_2$}
    \psline[linewidth=0.05pt]{->}(6.5,3)(6.5,2.5)

\end{pspicture*}}
\caption{$\Lambda_i,\ldots,\Lambda_n,\sigma(\Lambda_n),\ldots,\sigma(\Lambda_i)$}
\label{toughlemmafigure}
\end{figure}

One can see by considering cases that either $1)$
$\widetilde{\Lambda_{i-1}},\widetilde{\sigma(\Lambda_{i-1})}
\not\subset S$, $\widetilde{\Lambda_{i-1}} \not\subset
S^{\sigma}_{i-1}$, $\widetilde{\sigma(\Lambda_{i-1})} \not\subset
S_{i-1}$, or $2)$
$\widetilde{\Lambda_{i-1}},\widetilde{\sigma(\Lambda_{i-1})}
\not\subset S$, $\widetilde{\Lambda_{i-1}} \subset
S^{\sigma}_{i-1}$, $\widetilde{\sigma(\Lambda_{i-1})} \subset
S_{i-1}$, or $3)$ $\widetilde{\Lambda_{i-1}} \subset S$,
$\widetilde{\sigma(\Lambda_{i-1})} \not\subset S$, or $4)$
$\widetilde{\sigma(\Lambda_{i-1})} \subset S$,
$\widetilde{\Lambda_{i-1}} \not\subset S$.  One can also see by
considering cases that $i-1 = \min( \beta_1 , \beta_2 )$ (for
$\beta_1$,$\beta_2$ the two integers associated to some choice
edge $e$) if and only if we are in case $2$, $3$, or $4$. In
contrast to the $SL_3$ case, it is possible for neither
$\widetilde{\Lambda_{i-1}}$ nor
$\widetilde{\sigma(\Lambda_{i-1})}$ to be in $S$, while still
$\widetilde{\Lambda_{i-1}} \subset S_{i-1}^{\sigma}$ and
$\widetilde{\sigma(\Lambda_{i-1})} \subset S_{i-1}$.  To see this,
consider the case in which $p_r$ corresponds to an element of $W$
of length $3$ (pictured in Figure~\ref{toughlemmafigure}).  This
is the only situation in which case $2$ arises.

In case $1$ we can choose $\Lambda_{i-1}$ to be any alcove
adjacent to $\Lambda_i$ by $d_{i-1}$.  In this case, the dimension
of the space of choices for the construction increases by one, and
is therefore equal to $l(\Gamma^c_{(v_1,p_r)}) +
(n-(i-1))-1-n_{(F,i)}$.  We also have $n_{(F,i)} = n_{(F,i-1)}$.

Cases $3$ and $4$ only occur when $p_r$ corresponds to an element
of $W$ of length $2$.  We address case $3$.  The other case is
similar. We know $i-1=\min(\beta_1,\beta_2)$ for $\beta_1$ and
$\beta_2$ the two integers corresponding to some choice edge $e$.
If $F$ dictates a hard choice at $e$, we choose $\Lambda_{i-1}$ in
$S$. In this case there is no increase in the dimension of the
space of choices of the construction, and $n_{(F,i-1)} = n_{(F,i)}
+ 1$, so the dimension of the new space of choices is
$l(\Gamma^c_{(v_1,p_r)}) + (n-i)-1-n_{(F,i)} =
l(\Gamma^c_{(v_1,p_r)}) + (n-(i-1))-1-n_{(F,i-1)}$.  If $F$
dictates an easy choice at $e$, we choose $\Lambda_{i-1}$ to be
any alcove attached to $\Lambda_i$ at $d_{i-1}$, but not in $S$.
There is one dimension worth of such choices, so the dimension of
the space of choices of the construction increases by one.  We
have $n_{(F,i-1)} = n_{(F,i)}$, so the new dimension is
$l(\Gamma^c_{(v_1,p_r)}) + (n-i)-1-n_{(F,i)}+1 =
l(\Gamma^c_{(v_1,p_r)}) + (n-(i-1))-1-n_{(F,i-1)}$.

We now consider case $2$, which only occurs when $p_r$ corresponds
to an element of $W$ of length $3$. The construction
$\Lambda_i,\ldots,\Lambda_n,\sigma(\Lambda_n),\ldots,\sigma(\Lambda_i)$
is contained in an apartment, and has the general shape pictured
in Figure~\ref{toughlemmafigure}.  The dashed lines in this figure
represent the boundary of $S$.  Any choice of $\Lambda_{i-1}$
determines a $g(\Lambda_{i-1})$ attached to $\sigma(\Lambda_i)$
via $\sigma(d_{i-1})$, just by taking the alcove adjacent to
$\sigma(\Lambda_i)$ via $\sigma(d_{i-1})$ in the intersection of
all apartments containing $\Lambda_{i-1}$ and $S$.  We claim that
the number of choices of $\Lambda_{i-1}$ with $g(\Lambda_{i-1}) =
\sigma(\Lambda_{i-1})$ is non-zero and finite.  If this is true
and if $F$ requires a hard choice, take $\Lambda_{i-1}$ with
$g(\Lambda_{i-1}) = \sigma(\Lambda_{i-1})$.  Then dimension does
not increase, and is therefore equal to $l(\Gamma^c_{(v_1,p_r)}) +
(n-i)-1-n_{(F,i)}$, which is $l(\Gamma^c_{(v_1,p_r)}) +
(n-(i-1))-1-n_{(F,i-1)}$, since $n_{(F,i-1)} = n_{(F,i)} + 1$.  If
$F$ requires an easy choice, take $\Lambda_{i-1}$ with
$g(\Lambda_{i-1}) \neq \sigma(\Lambda_{i-1})$.  Then dimension
increases by one, and is therefore $l(\Gamma^c_{(v_1,p_r)}) +
(n-i)-1-n_{(F,i)}+1 = l(\Gamma^c_{(v_1,p_r)}) +
(n-(i-1))-1-n_{(F,i-1)}$, since $n_{(F,i-1)} = n_{(F,i)}$.

We now prove the claim of the previous paragraph.  We can identify
the set $\{ \Lambda_{i-1} \}$ with $\mathbb{A}^1$ over
$\bar{\mathbb{F}}_q$, where $\mathbb{F}_q$ is the residue field of
$F$.  We can identify the set $\{ \sigma(\Lambda_{i-1}) \}$ with
the set $\{ \Lambda_{i-1} \}$ (and therefore with $\mathbb{A}^1$)
using $g$.  The map $\sigma : \{ \Lambda_{i-1} \} \rightarrow \{
\sigma(\Lambda_{i-1}) \}$ given by the action of $\sigma$ on
$\mathcal{B}_{\infty}$ therefore gives a map $\psi : \mathbb{A}^1
\rightarrow \mathbb{A}^1$.  But $\sigma$ also acts on
$\mathbb{A}^1(\bar{\mathbb{F}}_q)$ as the (algebraic) Frobenius,
and one can show that if $\varphi : \mathbb{A}^1 \rightarrow
\mathbb{A}^1$ is defined by $\psi = \varphi \circ \sigma$, then
$\varphi$ is an algebraic isomorphism of $\mathbb{A}^1$.  So
$\varphi(x) = ax+b$ with $a \neq 0$.  The fixed points of $\psi$
correspond to $x \in \mathbb{A}^1$ such that $a \sigma(x) + b =
x$, which has exactly $q$ solutions since $\sigma(x) = x^q$.

\begin{figure}
\centerline{\input{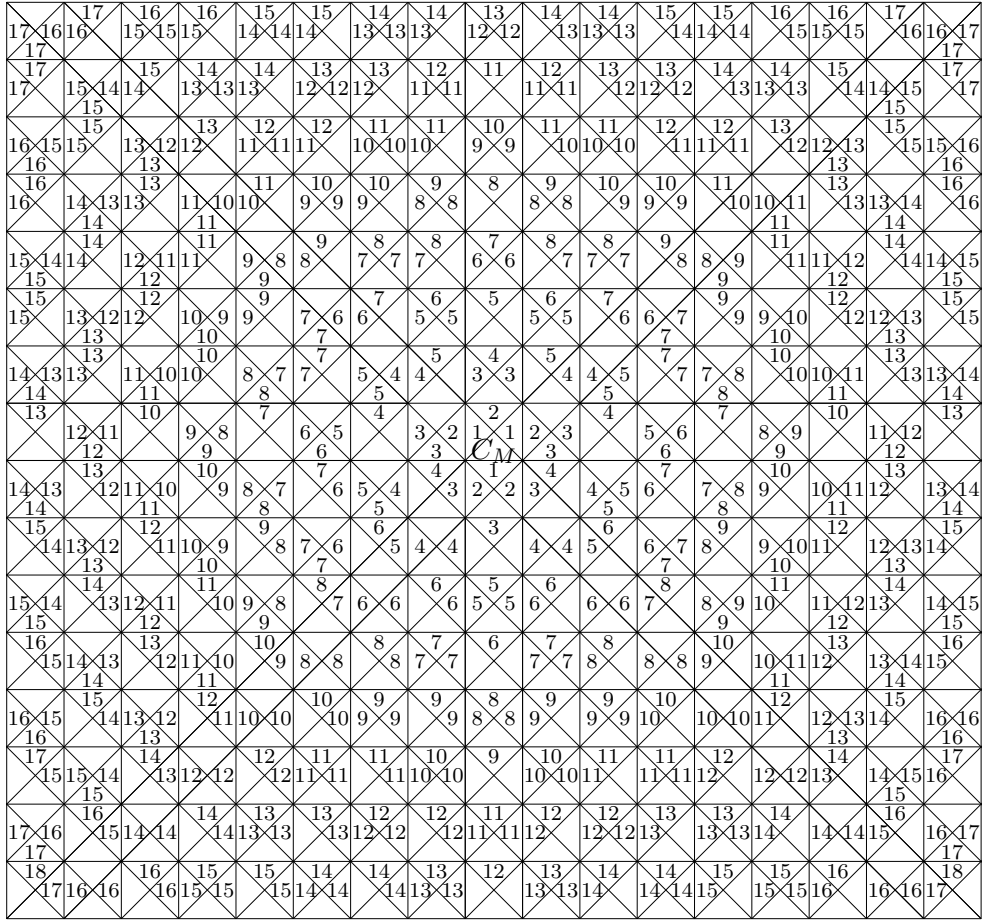}}
\caption{Main result in diagram
form for $Sp_4$} \label{Sp4result}
\end{figure}

So now we can compute emptiness/non-emptiness and dimension of the
$(v_1,v_2,p_r)$-piece of $X_{\tilde{w}}(\sigma)$ for each
$\tilde{w}$ by doing straightforward computations with
cf-dimension.  We did this for all $v_1$, and we reflected the
results across the line of symmetry of $C_M$.  We took maxima
whenever two numbers appeared in the same alcove.  The results of
this process can be found in Figure~\ref{Sp4result}.

In the course of the computation we observed that if the
$(v_1,v_2,p_r)$-piece of $X_{\tilde{w}}(\sigma)$ and the
$(v_1',v_2',p_r')$-piece of $X_{\tilde{w}}(\sigma)$ had different
dimensions, then exactly one of $v_1,v_1'$ was non-special, and
the corresponding piece had the smaller dimension.

\section{Application to
$\dim(X_{\tilde{w}}^K(\sigma))$}\label{dimK}

Let $\tilde{w} \in \tilde{W}$ and let $\mu$ be a dominant
cocharacter.  Let $\pi$ be the uniformizer in $F$.  The map
$G(L)/I \rightarrow G(L)/K$ gives a map $X_{\tilde{w}}^I(b\sigma)
\rightarrow X_{\mu(\pi)}^K(b\sigma)$ whenever $I\tilde{w}I \subset
K \mu (\pi) K$.  The non-empty fibers of this map are always
$K/I$, which has dimension equal to the length $\delta$ of the
longest element of the finite Weyl group, $W$. Further, any point
in $X_{\mu(\pi)}^K(b\sigma)$ is hit by a point in
$X_{\tilde{w}}^I(b\sigma)$ for some $\tilde{w}$ with $I \tilde{w}
I \subset K \mu(\pi) K$.  If $\mathcal{S}_{\mu(\pi)} \subset
\tilde{W}$ is defined so that $\coprod_{\tilde{w} \in
\mathcal{S}_{\mu(\pi)}} I \tilde{w} I = K \mu(\pi)K$, then we have
$\dim(X_{\mu(\pi)}^K(b\sigma))=\max_{\tilde{w} \in
\mathcal{S}_{\mu(\pi)}} (\dim(X_{\tilde{w}}^I(b\sigma))) -
\delta$.  We applied this formula to the cases $b=1$, $G=SL_2$,
$SL_3$, $Sp_4$, and found that $\dim(X_{\mu(\pi)}^K(\sigma)) =
\langle \mu , \rho \rangle$, where $\rho$ is half the sum of the
positive roots of $G$.  This result supports Rapoport's Conjecture
$5.10$ in \cite{R2}.

\section{A partial formula for $\dim(X_{\tilde{w}}^I(\sigma))$ for $SL_2$, $SL_3$ and $Sp_4$}
\label{formulasection}

Suppose $G$ is a simply-connected group and suppose $\tilde{w} \in
\tilde{W}$.  Let $\tilde{w} = tw$, where $w \in W$ and $t$ acts on
$A_M$ by translation.  Let $\eta_2(\tilde{w}) = \alpha \in W$,
where $\tilde{w}C_M$ is in the same Weyl chamber as $\alpha C_M$.
Let $\eta_1 : \tilde{W} \rightarrow W$ be the quotient map by the
subgroup of translations.  Let $S$ be the set of simple
reflections in $W$, and let $W_T$ be the subgroup of $W$ generated
by $T \subset S$.

Let $h_1,\ldots,h_{n+1}$ be the hyperplanes in $A_M$ that contain
one of the codimension-one sub-simplices of $C_M$.  Here $n$ is
the rank of $G$.  Let $h_i^{(j)}$ be the hyperplanes in $A_M$
parallel to $h_i$, with $h_i^{(0)} = h_i$.  Choose $h_i^{(1)}$ to
be as close as possible to $h_i$, but on the other side of $C_M$.
We define the {\em union of shrunken Weyl chambers} to be the set
of all alcoves that are not between $h_i^{(0)}$ and $h_i^{(1)}$
for any $i$.

If $\tilde{w}C_M$ is in the union of shrunken Weyl chambers and if
$G=SL_2$, $SL_3$, or $Sp_4$, then $X_{\tilde{w}}(\sigma)$ is
non-empty if and only if $\eta_2(\tilde{w})^{-1} \eta_1(\tilde{w})
\eta_2(\tilde{w}) \in W \setminus \cup_{T \subset S} W_T$, and in
this case
$$\dim(X_{\tilde{w}}(\sigma)) = \frac{l_{\tilde{W}}(\tilde{w}) +
l_{W}(\eta_2(\tilde{w})^{-1} \eta_1(\tilde{w}) \eta_2(\tilde{w})
)}{2}.$$ Here $l_W$ is length in $W$ and $l_{\tilde{W}}$ is length
in $\tilde{W}$, as Coxeter groups.

One can examine Figures~\ref{SL3result} and~\ref{Sp4result} to see
that the above statement holds for $SL_3$ and $Sp_4$.  It also
holds for $SL_2$. Note, though, that the new formula says nothing
about the dimension or emptiness/non-emptiness of
$X_{\tilde{w}}(\sigma)$ for $\tilde{w}$ not in the union of the
shrunken Weyl chambers. Figures~\ref{SL3result}
and~\ref{Sp4result} give this information for $SL_3$ and $Sp_4$.
The complement of the union of the shrunken Weyl chambers for
$SL_2$ is just $C_M$, an easy special case.

The above formula might not hold for $SL_4$.  One problem is that,
for $SL_2$, $SL_3$, and $Sp_4$, there are other ways to specify
the set $W \setminus \cup_{T \subset S} W_T$. In particular, $W
\setminus \cup_{T \subset S} W_T = \{ w \in W : l_W(w) \geq
\rank(G) \}$ for these groups, and the $SL_4$ analogues of these
two sets are not the same.  We think the first formulation is more
likely to be appropriate for a general statement.

\section{Related results}\label{RelatedResults}

Some of the emptiness/non-emptiness results of this paper were
also obtained using other methods in the author's Ph.D. thesis
\cite{Re1}. These other methods are more computationally
intensive, and do not provide dimension information, but they
extend to some extent to $b \neq 1$.  Some of the results from
\cite{Re1} can be combined to suggest the below conjecture.

We restrict $G$ to be one of the groups $SL_2$, $SL_3$, and
$Sp_4$.  Let $D$ be a Weyl chamber in $A_M$, and let $D'$ be the
intersection of $D$ with the union of the shrunken Weyl chambers.
Then we call $D'$ a shrunken Weyl chamber.  Let $b$ be a
representative of a $\sigma$-conjugacy class that meets the main
torus of $G$. We can choose $b$ so that it acts on $A_M$ by
translation, and such that $bC_M$ is in the main Weyl chamber. We
define the {\em $b$-shifted shrunken Weyl chamber} associated to
$D$ to be $wbw^{-1}D'$, where $w \in W$ is the Weyl group element
corresponding to the Weyl chamber $D$.

\begin{conj}
If $b$ and $G$ are restricted as above, and if $\tilde{w}C_M$ is
in the union of the $b$-shifted shrunken Weyl chambers, then
$X_{\tilde{w}}(b \sigma)$ is non-empty if and only if
$\eta_2(\tilde{w})^{-1} \eta_1(\tilde{w}) \eta_2(\tilde{w}) \in W
\setminus \cup_{T \subset S} W_T$.
\end{conj}

This conjecture is shown to hold true in \cite{Re1} for several
values of $b$.  Information about $\tilde{w}$ not in the
$b$-shifted shrunken Weyl chambers is also given in \cite{Re1} for
the same $b$ values, but we have been unable to describe these
results with a formula.

\section{Acknowledgements}

We thank Robert Kottwitz for suggesting the problem. We also thank
Michael Rapoport for helpful suggestions and a helpful
conversation, and Joel Cohen for his support and understanding
during the writing process.  We acknowledge with gratitude U.S.
National Science Foundation grant DEB 9981552 for supporting some
of the research.

\end{document}